\pdfoutput=1
\documentclass{amsart}
\usepackage{ascmac}
\usepackage{amsmath}
\usepackage{amscd}
\usepackage{amssymb}
\usepackage{hyperref}
\usepackage[dvipdfmx]{graphicx}
\usepackage{here}
\usepackage{color}
\usepackage{hyperref}

\theoremstyle{plain}
\newtheorem{lemma}{Lemma}

\newtheorem{theorem}{Theorem}

\newtheorem{fact}{Fact}
\newcommand{\LHS}{\operatorname{LHS}}
\newcommand{\RHS}{\operatorname{RHS}}
\newcommand{\sgn}{\operatorname{sgn}}
\newcommand{\BMb}{
\begin{picture}(0,0)
\put(0,11){$\epsilon$}
\qbezier(-10,0)(0,10)(-10,20)
\qbezier(15,0)(5,10)(15,20)
\put(-5,10){\vector(4,-1){15}}
\put(10,10){\circle*{2}}
\end{picture}
}
\newcommand{\BMc}{
\begin{picture}(0,0)
\put(-3,8){$\to$}
\end{picture}}
\newcommand{\BMa}{
  \begin{picture}(0,0)
\put(0,14){$\epsilon$}
\qbezier(-10,0)(0,10)(-10,20)
\qbezier(15,0)(5,10)(15,20)
\put(-5,10){\vector(4,1){15}}
\put(10,10){\circle*{2}}
\end{picture}
}

\newcommand{\bMb}{
\begin{picture}(0,0)
 \put(0,10){$\eta$}
\qbezier(-10,0)(0,10)(-10,20)
\qbezier(15,0)(5,10)(15,20)
\put(11,5){\vector(-4,1){16}}
\put(10,10){\circle*{2}}
\end{picture}}
\newcommand{\bMa}{
\begin{picture}(0,0)
 \put(0,15){$\eta$}
\qbezier(-10,0)(0,10)(-10,20)
\qbezier(15,0)(5,10)(15,20)
\put(11,14){\vector(-4,-1){16}}
\put(10,10){\circle*{2}}
\end{picture}
}
\newcommand{\sign}{\operatorname{sign}}
\newcommand{\RR}{
\begin{picture}(50,8)
\put(3,0.5){\tiny$i$}
\put(23,1.5){\tiny$j$}
\put(43,0.5){\tiny$k$}
\put(5,8){\circle*{2}}
\put(25,8){\circle*{2}}
\put(45,8){\circle*{2}}
\put(5,3){\circle{10}}
\put(10,3){\vector(1,0){10}}
\put(25,3){\circle{10}}
\put(30,3){\vector(1,0){10}}
\put(45,3){\circle{10}}
\end{picture}
}
\newcommand{\RRsign}{
\begin{picture}(50,8)
\put(3,0.5){\tiny$i$}
\put(23,1.5){\tiny$j$}
\put(43,0.5){\tiny$k$}
\put(5,8){\circle*{2}}
\put(25,8){\circle*{2}}
\put(45,8){\circle*{2}}
\put(5,3){\circle{10}}
\put(10,3){\vector(1,0){10}}
\put(10,7){\tiny$\epsilon$}
\put(30,7){\tiny$\eta$}
\put(25,3){\circle{10}}
\put(30,3){\vector(1,0){10}}
\put(45,3){\circle{10}}
\end{picture}
}

\newcommand{\LR}{
\begin{picture}(50,8)
\put(5,8){\circle*{2}}
\put(25,8){\circle*{2}}
\put(45,8){\circle*{2}}
\put(5,3){\circle{10}}
\put(3,0.5){\tiny$i$}
\put(23,1.5){\tiny$j$}
\put(43,0.5){\tiny$k$}
\put(20,3){\vector(-1,0){10}}
\put(25,3){\circle{10}}
\put(30,3){\vector(1,0){10}}
\put(45,3){\circle{10}}
\end{picture}
}
\newcommand{\LRikj}{
\begin{picture}(50,8)
\put(5,8){\circle*{2}}
\put(25,8){\circle*{2}}
\put(45,8){\circle*{2}}
\put(5,3){\circle{10}}
\put(3,0.5){\tiny$i$}
\put(23,1.5){\tiny$k$}
\put(43,0.5){\tiny$j$}
\put(20,3){\vector(-1,0){10}}
\put(25,3){\circle{10}}
\put(30,3){\vector(1,0){10}}
\put(45,3){\circle{10}}
\end{picture}
}
\newcommand{\LRikjsign}{
\begin{picture}(50,8)
\put(13,7){\tiny$\epsilon$}
\put(33,7){\tiny$\eta$}
\put(5,8){\circle*{2}}
\put(25,8){\circle*{2}}
\put(45,8){\circle*{2}}
\put(5,3){\circle{10}}
\put(3,0.5){\tiny$i$}
\put(23,1.5){\tiny$k$}
\put(43,0.5){\tiny$j$}
\put(20,3){\vector(-1,0){10}}
\put(25,3){\circle{10}}
\put(30,3){\vector(1,0){10}}
\put(45,3){\circle{10}}
\end{picture}
}

\newcommand{\RL}{
\begin{picture}(50,8)
\put(3,0.5){\tiny$i$}
\put(23,1.5){\tiny$j$}
\put(43,0.5){\tiny$k$}
\put(5,8){\circle*{2}}
\put(25,8){\circle*{2}}
\put(45,8){\circle*{2}}
\put(5,3){\circle{10}}
\put(10,3){\vector(1,0){10}}
\put(25,3){\circle{10}}
\put(40,3){\vector(-1,0){10}}
\put(45,3){\circle{10}}
\end{picture}
}
\newcommand{\RLsign}{
\begin{picture}(50,8)
\put(13,7){\tiny$\epsilon$}
\put(33,7){\tiny$\eta$}
\put(3,0.5){\tiny$i$}
\put(23,1.5){\tiny$j$}
\put(43,0.5){\tiny$k$}
\put(5,8){\circle*{2}}
\put(25,8){\circle*{2}}
\put(45,8){\circle*{2}}
\put(5,3){\circle{10}}
\put(10,3){\vector(1,0){10}}
\put(25,3){\circle{10}}
\put(40,3){\vector(-1,0){10}}
\put(45,3){\circle{10}}
\end{picture}
}

\newcommand{\RLkji}{
\begin{picture}(50,8)
\put(3,0.5){\tiny$k$}
\put(23,1.5){\tiny$j$}
\put(43,0.5){\tiny$i$}
\put(5,8){\circle*{2}}
\put(25,8){\circle*{2}}
\put(45,8){\circle*{2}}
\put(5,3){\circle{10}}
\put(10,3){\vector(1,0){10}}
\put(25,3){\circle{10}}
\put(40,3){\vector(-1,0){10}}
\put(45,3){\circle{10}}
\end{picture}
}
\newcommand{\RLkjisign}{
\begin{picture}(50,8)
\put(13,7){\tiny$\eta$}
\put(33,7){\tiny$\epsilon$}
\put(3,0.5){\tiny$k$}
\put(23,1.5){\tiny$j$}
\put(43,0.5){\tiny$i$}
\put(5,8){\circle*{2}}
\put(25,8){\circle*{2}}
\put(45,8){\circle*{2}}
\put(5,3){\circle{10}}
\put(10,3){\vector(1,0){10}}
\put(25,3){\circle{10}}
\put(40,3){\vector(-1,0){10}}
\put(45,3){\circle{10}}
\end{picture}
}

\newcommand{\LL}{
\begin{picture}(50,8)
\put(3,0.5){\tiny$i$}
\put(23,1.5){\tiny$j$}
\put(43,0.5){\tiny$k$}
\put(5,8){\circle*{2}}
\put(25,8){\circle*{2}}
\put(45,8){\circle*{2}}
\put(5,3){\circle{10}}
\put(20,3){\vector(-1,0){10}}
\put(25,3){\circle{10}}
\put(40,3){\vector(-1,0){10}}
\put(45,3){\circle{10}}
\end{picture}
}

\newcommand{\LLkji}{
\begin{picture}(50,8)
\put(3,0.5){\tiny$k$}
\put(23,1.5){\tiny$j$}
\put(43,0.5){\tiny$i$}
\put(5,8){\circle*{2}}
\put(25,8){\circle*{2}}
\put(45,8){\circle*{2}}
\put(5,3){\circle{10}}
\put(20,3){\vector(-1,0){10}}
\put(25,3){\circle{10}}
\put(40,3){\vector(-1,0){10}}
\put(45,3){\circle{10}}
\end{picture}
}
\newcommand{\LLkjisign}{
\begin{picture}(50,8)
\put(3,0.5){\tiny$k$}
\put(23,1.5){\tiny$j$}
\put(43,0.5){\tiny$i$}
\put(5,8){\circle*{2}}
\put(25,8){\circle*{2}}
\put(45,8){\circle*{2}}
\put(5,3){\circle{10}}
\put(13,5){\tiny$\eta$}
\put(20,3){\vector(-1,0){10}}
\put(25,3){\circle{10}}
\put(33,5){\tiny$\epsilon$}
\put(40,3){\vector(-1,0){10}}
\put(45,3){\circle{10}}
\end{picture}
}

\newcommand{\LRjki}{
\begin{picture}(50,8)
\put(5,8){\circle*{2}}
\put(25,8){\circle*{2}}
\put(45,8){\circle*{2}}
\put(5,3){\circle{10}}
\put(3,1.5){\tiny$j$}
\put(23,1.5){\tiny$k$}
\put(43,1.5){\tiny$i$}
\put(20,3){\vector(-1,0){10}}
\put(25,3){\circle{10}}
\put(30,3){\vector(1,0){10}}
\put(45,3){\circle{10}}
\end{picture}
}
\newcommand{\LRjkisign}{
\begin{picture}(50,8)
\put(13,7){\tiny$\eta$}
\put(33,7){\tiny$\epsilon$}
\put(5,8){\circle*{2}}
\put(25,8){\circle*{2}}
\put(45,8){\circle*{2}}
\put(5,3){\circle{10}}
\put(3,1.5){\tiny$j$}
\put(23,1.5){\tiny$k$}
\put(43,1.5){\tiny$i$}
\put(20,3){\vector(-1,0){10}}
\put(25,3){\circle{10}}
\put(30,3){\vector(1,0){10}}
\put(45,3){\circle{10}}
\end{picture}
}

\newcommand{\RRenr}{ 
\begin{picture}(50,20)(0,-4)
\put(3,0.5){\tiny$i$}
\put(23,1.5){\tiny$j$}
\put(43,0.5){\tiny$k$}
\put(5,8){\circle*{2}}
\put(25,8){\circle*{2}}
\put(45,8){\circle*{2}}
\put(5,3){\circle{10}}
\put(13,-1){\tiny$+$}
\qbezier(7,-2)(14,-6)(23,-2) 
\put(24,-2){\vector(3,1){1}} 
\put(25,3){\circle{10}}
\put(33,-1){\tiny$+$}
\qbezier(27,-2)(36,-6)(43,-2) 
\put(44,-2){\vector(3,1){1}} 
\put(45,3){\circle{10}}
\put(23,-6.5){\tiny$-$}
\qbezier(3,-2)(30,-14)(48,-1.5) 
\put(47.5,-2){\vector(1,1){1}} 
\end{picture}
}

\newcommand{\RRenrS}{ 
\begin{picture}(50,25)(0,-8)
\put(3,0.5){\tiny$i$}
\put(23,1.5){\tiny$j$}
\put(43,0.5){\tiny$k$}
\put(5,8){\circle*{2}}
\put(25,8){\circle*{2}}
\put(45,8){\circle*{2}}
\put(5,3){\circle{10}}
\put(14,-4){\tiny$+$}
\qbezier(3,-2)(19,-11)(28,-1.5) 
\put(28,-1.5){\vector(2,3){1}} 
\put(25,3){\circle{10}}
\put(32,-4){\tiny$+$}
\qbezier(22,-1.5)(33,-11)(47,-2) 
\put(47,-2){\vector(1,1){1}} 
\put(45,3){\circle{10}}
\put(23,-9){\tiny$-$}
\qbezier(7,-2)(27,-18)(44,-1.75) 
\put(43.5,-2.5){\vector(1,1){1}} 
\end{picture}
}
\newcommand{\RRen}{ 
\begin{picture}(50,20)(0,-4)
\put(3,0.5){\tiny$i$}
\put(23,1.5){\tiny$j$}
\put(43,0.5){\tiny$k$}
\put(5,8){\circle*{2}}
\put(25,8){\circle*{2}}
\put(45,8){\circle*{2}}
\put(5,3){\circle{10}}
\put(13,-1){\tiny$+$}
\qbezier(7,-2)(14,-6)(23,-2) 
\put(24,-2){\vector(3,1){1}} 
\put(25,3){\circle{10}}
\put(33,-1){\tiny$+$}
\qbezier(27,-2)(36,-6)(43,-2) 
\put(44,-2){\vector(3,1){1}} 
\put(45,3){\circle{10}}
\end{picture}
}

\newcommand{\RRer}{ 
\begin{picture}(50,20)(0,-4)
\put(3,0.5){\tiny$i$}
\put(23,1.5){\tiny$j$}
\put(43,0.5){\tiny$k$}
\put(5,8){\circle*{2}}
\put(25,8){\circle*{2}}
\put(45,8){\circle*{2}}
\put(5,3){\circle{10}}
\put(13,-1){\tiny$+$}
\qbezier(7,-2)(14,-6)(23,-2) 
\put(24,-2){\vector(3,1){1}} 
\put(25,3){\circle{10}}
\put(45,3){\circle{10}}
\put(23,-6.5){\tiny$-$}
\qbezier(3,-2)(30,-14)(48,-1.5) 
\put(47.5,-2){\vector(1,1){1}} 
\end{picture}
}

\newcommand{\RRnr}{ 
\begin{picture}(50,20)(0,-4)
\put(3,0.5){\tiny$i$}
\put(23,1.5){\tiny$j$}
\put(43,0.5){\tiny$k$}
\put(5,8){\circle*{2}}
\put(25,8){\circle*{2}}
\put(45,8){\circle*{2}}
\put(5,3){\circle{10}}
\put(25,3){\circle{10}}
\put(33,-1){\tiny$+$}
\qbezier(27,-2)(36,-6)(43,-2) 
\put(44,-2){\vector(3,1){1}} 
\put(45,3){\circle{10}}
\put(23,-6.5){\tiny$-$}
\qbezier(3,-2)(30,-14)(48,-1.5) 
\put(47.5,-2){\vector(1,1){1}} 
\end{picture}
}

\newcommand{\RRenS}{ 
\begin{picture}(50,25)(0,-8)
\put(3,0.5){\tiny$i$}
\put(23,1.5){\tiny$j$}
\put(43,0.5){\tiny$k$}
\put(5,8){\circle*{2}}
\put(25,8){\circle*{2}}
\put(45,8){\circle*{2}}
\put(5,3){\circle{10}}
\put(14,-4){\tiny$+$}
\qbezier(3,-2)(19,-11)(28,-1.5) 
\put(28,-1.5){\vector(2,3){1}} 
\put(25,3){\circle{10}}
\put(32,-4){\tiny$+$}
\qbezier(22,-1.5)(33,-11)(47,-2) 
\put(47,-2){\vector(1,1){1}} 
\put(45,3){\circle{10}}
\end{picture}
}

\newcommand{\RRerS}{ 
\begin{picture}(50,25)(0,-8)
\put(3,0.5){\tiny$i$}
\put(23,1.5){\tiny$j$}
\put(43,0.5){\tiny$k$}
\put(5,8){\circle*{2}}
\put(25,8){\circle*{2}}
\put(45,8){\circle*{2}}
\put(5,3){\circle{10}}
\put(14,-4){\tiny$+$}
\qbezier(3,-2)(19,-11)(28,-1.5) 
\put(28,-1.5){\vector(2,3){1}} 
\put(25,3){\circle{10}}
\put(45,3){\circle{10}}
\put(23,-9){\tiny$-$}
\qbezier(7,-2)(27,-18)(44,-1.75) 
\put(43.5,-2.5){\vector(1,1){1}} 
\end{picture}
}

\newcommand{\RRnrS}{ 
\begin{picture}(50,25)(0,-8)
\put(3,0.5){\tiny$i$}
\put(23,1.5){\tiny$j$}
\put(43,0.5){\tiny$k$}
\put(5,8){\circle*{2}}
\put(25,8){\circle*{2}}
\put(45,8){\circle*{2}}
\put(5,3){\circle{10}}
\put(25,3){\circle{10}}
\put(32,-4){\tiny$+$}
\qbezier(22,-1.5)(33,-11)(47,-2) 
\put(47,-2){\vector(1,1){1}} 
\put(45,3){\circle{10}}
\put(23,-9){\tiny$-$}
\qbezier(7,-2)(27,-18)(44,-1.75) 
\put(43.5,-2.5){\vector(1,1){1}} 
\end{picture}
}

\newcommand{\sub}{\operatorname{Sub}}
\newcommand{\SGL}{\langle S,G_{l} \rangle}
\newcommand{\SGLs}{\langle S, G_{r} \rangle}
\newcommand{\SZ}{\langle S, z^{(l)}_i \rangle}
\newcommand{\SZs}{\langle S, z^{(r)}_{i} \rangle}
\newcommand{\Sumi}{\sum_{i=0}^{3}}

\newcommand{\SumZGD}{\sum_{z_i \in \sub(G_D)}}
\newcommand{\SumZGDs}{\sum_{z'_i \in \sub(G_D')}}

\title{Variations of Milnor's triple linking number}
\author{Kamolphat Intawong and Noboru Ito}
\begin{document}
\begin{abstract}
Topological polymers have various topological types, and they are expressed by graphs.  However, the Jones polynomial, we have a difficulty to compute it; computational time is growing exponentially with respect to the crossing number.   
The simplest Vassiliev invariant is the linking number and thus we will seek a next simple one is as the Milnor's triple linking number.  In this paper, we introduce simple Gauss diagram formulas of Vassiliev invariants of Milnor type.   
These are non-torsion valued, whereas the base-point-free Milnor's triple linking number is usually torsion-valued.   
\end{abstract}
\thanks{MSC2020:57K10; 57K16}
\keywords{Gauss diagram; Milnor invariant; links; Vassiliev invariant}
\date{May 9, 2022}
\maketitle
\section{Introduction}\label{sec:intro}
Polymers have various topological types and they are expressed by graphs.   
\emph{Topological polymers} indicate experimental / theoretical meanings.  One of them is a polymer of complex chemical connectivity that is expressed with graphs and synthesized in experiments; the other is a polymer with nontrivial topology realized as embedded spatial graphs in the three dimensional space.   

Nowadays since chemists synthesized polymers whose are some topologically distinct graphs, mathematicians are requested functions to detect their topological complexity with computational efficiency.  
As is well known, in practical application of the Jones polynomial of links, we have a difficulty to compute it; computational time is growing exponentially with respect to the crossing number of link diagrams.  

Deguchi \cite{Deguchi1994} gave an algorithm by the Taylor expansion at $q=1$ for the Jones polynomial $V_K (q)$: 
\[
V_K (q) = 1 + v_2 (K) \epsilon^2 + v_3 (K) \epsilon^3 + \dots , 
\]
where the coefficient $v_n (K)$ is called a Vassiliev invariant of degree $n$ and is expressed by a Gauss diagram formula in theory.  It is considered that Gauss diagram formulas are the simplest for the computations purpose.  Moreover, the simplest Vassiliev invariant is the linking number and a next simple one is the Milnor's triple linking number.  Thus, we search simple Gauss diagram formulas of Vassiliev invariants of Milnor type.   

In this paper, we find new Gauss diagram formulas that are non-torsion triple linking number that is similar to  Milnor-type, which are base-point-free and link homotopy invariants.  The purpose 
is of two things.  One of them is to give a sufficient simple  for the computation of applied mathematics; the other is to give  patterns of the affine index polynomial derived from a triple product \cite{Ito2022}.   

\section{Main Result}\label{sec:main}
{The notations of Gauss diagrams and Gauss diagrams obey \cite{Ostlund2004} }
Theorem~\ref{main} is the first result of this paper.   

\newcommand{\Psigma}{a_\sigma \RR+b_\sigma \LL +c_\sigma \RL +d_\sigma\ \LR}
\newcommand{\SPsigma}{\sum_{\sigma \in S_3}\sgn(\sigma)P^\sigma}
\newcommand{\F}{ f(a_\sigma,b_\sigma,c_\sigma,d_\sigma)}
\begin{theorem}\label{main}
{\color{black}{Let $L$ be an ordered three-component link, $G_L$ a Gauss diagram of $L$, $\sigma$ a permutation 
$\left(\begin{matrix}
1& 2 & 3 \\ i & j & k
\end{matrix}\right)
$ of the circles.  Let  
\[
P^\sigma=\Psigma,
\]
\[
S=\SPsigma,~{\textrm{and}}~ \F(L)=\langle S,G_L \rangle.  
\]
Suppose one of the following conditions: 
\begin{enumerate}
\item\label{I} $a_{ijk}=b_{kji}$, others are 0  for any $\sigma=\left(\begin{smallmatrix}
1& 2 & 3 \\ i & j & k
\end{smallmatrix}\right)
$, i.e. for a constant $\lambda=a_{ijk}$,  
\[
\F(L)= \lambda \left \langle \sum_{\sigma} \sgn(\sigma)  \left(\RR + \LLkji \right), G_L \right \rangle \quad.  
\]
\item\label{J} $c_{ijk}=c_{kji} = -d_{jki}= -d_{ikj}$, others are 0  for any even permutation $\sigma=\left(\begin{smallmatrix}
1& 2 & 3 \\ i & j & k
\end{smallmatrix}\right)
$, i.e. for a constant $\lambda=c_{ijk}$,  
\[
\F(L)= \lambda \left\langle \RL - \LRjki + \RLkji - \LRikj, G_L \right\rangle.
\]
\end{enumerate}
Then, $f(a_\sigma,b_\sigma,c_\sigma,d_\sigma)(L)$ is invariant under Reidemeister moves $($isotopy, strictly speaking link homotopy$)$. 
}}
\end{theorem}
\begin{fact}[\cite{PolyakViro1994}, \cite{Ostlund2004}]\label{Milnor}
Let $L$ be an ordered $3$-component link which components are $k_i$ $(i=1, 2, 3)$.  Then  
$\frac{1}{6}f(2, 2, 1, 1)(L) \mod \gcd(lk(k_2, k_3), lk(k_1, k_3), lk(k_1, k_2))$ is the Milnor link invariant of $L$.   
\end{fact}
\begin{theorem}\label{IndeTHM}
The family (\ref{I}) of invariants is  independent of the Milnor's $\mu_{123}$.  
\end{theorem}
The above results concludes that the generating function $\F$ includes $\mu_{123}$ and different invariants from $\mu_{123}$.  This fact implies that $\F$ is stronger than $\mu_{123}$ in a sense.  
\section{Proof of Theorem~\ref{main}}
In this section, we freely use \"{O}stlund's terminologies \cite{Ostlund2004} especially including Gauss diagrams, arrow diagrams, Gauss/arrow diagram fragments, and their Reidemeister moves  \cite[Section~1.6, Section~4, and Table~1]{Ostlund2004} except for replacing $\Omega_{+---}$ in \cite[Table~1]{Ostlund2004} with $\Omega_{+-+-}$ in \cite[Figure~9]{Ito2022} \footnote{This notation  obeys \cite{Ostlund2001PhD}.  
The authors have not proved that \cite[Table~1]{Ostlund2004} is really a generating set but have found the  generating set after applying the replacement of $\Omega_{+---}$ in \cite[Table~1]{Ostlund2004} with $\Omega_{+-+-}$ in \cite[Figure~9]{Ito2022}.}.  

In the rest of this section, without loss of generality, we suppose that $\lambda=1$.  
It is clear that each function as in the statement satisfy the invariance under the first Reidemeister move.
\subsection{Invariances under the Reidemeister moves in one component}
It is clear that our functions are invariant under any Reidemeister moves in one component.  
\subsection{Invariances under the second and third Reidemeister moves in two components}
The second (third,~resp) Reidemeister invariance in two components is given by \cite[Section~4.5 (2)]{Ostlund2004} (\cite[Section~4.5 (3)]{Ostlund2004},~resp.).  
  
\subsection{Invariance under $\Omega_{III+-+3}$}
In this section, we use word-theoretic notation of \cite[Notation~4]{ItoOyamaguchi2021}.    
As a preliminary of the proof, we recall two lemmas as analogues of  \cite[Lemmas~3 and 4]{ItoOyamaguchi2021}.  These are given by \cite[Section~4.5]{Ostlund2004}.  
\begin{lemma}\label{ZeroArrow}
For any linear sum $A$ of arrow diagrams, 
\[
\sum_{z^{(l)}_0 \in \sub(D^{\Omega_{III+-+3}}_l)} \sign (z^{(l)}_0) (A, z^{(l)}_0) = \sum_{z_0^{(r)}  \in \sub(D^{\Omega_{III+-+3}}_r)} \sign (z_0^{(r)}) (A, z_0^{(r)}) \quad.  
\]
\end{lemma}
\begin{lemma}\label{OneArrow}
For any linear sum $A$ of arrow diagrams, 
\[
\sum_{z^{(l)}_1 \in \sub(D^{\Omega_{III+-+3}}_l)} \sign (z^{(l)}_1) (A, z^{(l)}_1) = \sum_{z_1^{(r)}  \in \sub(D^{\Omega_{III+-+3}}_r)} \sign (z_1^{(r)}) (A, z_1^{(r)}) \quad.  
\]
\end{lemma}

\begin{itemize}
\item Case~(\ref{I}).  
\begin{align*}
&\SGLs - \SGL =  \Sumi \left(\SumZGD \SZs - \SumZGDs \SZ \right)\\
&\stackrel{\textrm{Lemmas~\ref{ZeroArrow}, \ref{OneArrow}}}{=} \sum_{z_2 \in \sub(D_r)} \langle S, z^{(r)}_2 \rangle - \sum_{z^{(l)}_2 \in \sub(D_l)} \langle S, z^{(l)}_{2} \rangle\\
&= \sum_{z_2 \in \sub(D_r)} \langle S, z_{2a} + z_{2b} + z_{2c} - z'_{2a} - z'_{2b} - z'_{2c}  \rangle \\
&= \sum_{z_2 \in \sub(D_r)} \sum_{\sigma} \sgn(\sigma) \langle P^{\sigma}, z_{2a} + z_{2b} + z_{2c} - z'_{2a} - z'_{2b} - z'_{2c}  \rangle \\
&= \sum_{z_2 \in \sub(D_r)} \sum_{\sigma} \sgn(\sigma) \langle \RR + \LLkji, z_{2a} + z_{2b} + z_{2c} - z'_{2a} - z'_{2b} - z'_{2c}  \rangle \\
&= 0,  
\end{align*}
where 
$z'_{2b}$ $=$ $\RRen$ , $z'_{2a}$ $=$ $\RRer$, $z'_{2c}$ $=$ $\RRnr$ , $z_{2b}$ $=$ $\RRenS$, \\
$z_{2a}$ $=$ $\RRerS$, $z_{2c}$ $=$ $\RRnrS$ (whose indices are as in \cite[Proof of Theorem~1]{Ito2019}) \\
which are the 2-{\color{black}{arrow}} sub-Gauss diagrams of a 3-component link before and after applying {\color{black}{$\Omega_{III+-+3}$}}. Spefically, $\RRenr$ and $\RRenrS$. 

\begin{align*}
&\sum_{z_2 \in \sub(D_r)} \sum_{\sigma} \sgn(\sigma) \langle \RR + \LLkji, z_{2a} + z_{2b} + z_{2c} - z'_{2a} - z'_{2b} - z'_{2c}  \rangle \\
&= \sum_{z_2 \in \sub(D_r)} \sum_{\sigma} \sgn(\sigma) \langle \RR + \LLkji, z_{2b}  - z'_{2b}   \rangle \\
&=0.
\end{align*}

\item Case~(\ref{J}). 
\\ The same argument as that of Case~(\ref{I}) is applied to Case~(\ref{J}). Then,
\begin{align*}
&\sum_{z_2 \in \sub(D_r)}  \langle \RL - \LRjki + \RLkji - \LRikj, \\
&\qquad\qquad\qquad z_{2a} + z_{2b} + z_{2c} - z'_{2a} - z'_{2b} - z'_{2c}  \rangle \\
&=\sum_{z_2 \in \sub(D_r)}  \langle \RL - \LRjki, - z'_{2a} - z'_{2c}  \rangle \\
&+\sum_{z_2 \in \sub(D_r)}  \langle \RLkji - \LRikj, z_{2a} + z_{2c}  
  \rangle \\
&=0.
\end{align*}
\end{itemize}
\subsection{Invariance under base point moves}
It is sufficient to check the invariances under 
\[\]
\qquad\qquad $
\BMb \qquad \BMc \qquad \BMa
$ \qquad (Case~A) and \qquad  
$
\bMb \qquad \BMc \qquad \bMa
$ \qquad (Case~B).  
For each case, the diagram of the left (right,~resp.) hand side is denoted by $D_l$ ($D_r$,~resp.).  In the following, we see comparisons of   $D_l$ and $D_r$ after applying Lemmas~\ref{ZeroArrow} and \ref{OneArrow}.   
\subsection{Case~A} 
\begin{itemize}
\item Case~(\ref{I}).  
\begin{align*}
\LHS - \RHS =&\sum_{z_2 \in \sub(D_l)} \sum_{\sigma} \sgn(\sigma) \langle \RR + \LLkji,  \RRsign  \rangle\\
&- \sum_{z_2 \in \sub(D_r)} \sum_{\sigma} \sgn(\sigma) \langle \RR + \LLkji,  \LLkjisign \rangle \\
&=0.
\end{align*}
\item Case~(\ref{J}). 
\begin{align*}
&\LHS - \RHS\\ 
=&\sum_{z_2 \in \sub(D_l)}  \langle \RL - \LRjki + \RLkji - \LRikj, \RLsign   \rangle\\
&- \sum_{z_2 \in \sub(D_r)} \langle \RL - \LRjki + \RLkji - \LRikj,   \RLkjisign \rangle \\
&=0.
\end{align*}
\end{itemize}
\subsection{Case~B}
\begin{itemize}
\item Case~(\ref{I}). 
\begin{align*}
\LHS - \RHS =&\sum_{z_2 \in \sub(D_l)} \sum_{\sigma} \sgn(\sigma) \langle \RR + \LLkji,  \LLkjisign  \rangle\\
&- \sum_{z_2 \in \sub(D_r)} \sum_{\sigma} \sgn(\sigma) \langle \RR + \LLkji, \RRsign \rangle \\
&=0.
\end{align*}
\item Case~(\ref{J}). 
\begin{align*}
&\LHS - \RHS\\ 
=&\sum_{z_2 \in \sub(D_l)}  \langle \RL - \LRjki + \RLkji - \LRikj, \LRjkisign   \rangle\\
&- \sum_{z_2 \in \sub(D_r)} \langle \RL - \LRjki + \RLkji - \LRikj,   \LRikjsign \rangle \\
&=0.
\end{align*}
\end{itemize}

\section{Proof of Theorem~\ref{IndeTHM}}\label{sec:ex}
In this section, we will show that $\F$ is independent of the  Milnor's $\mu_{123}$. More specifically, we claim that there is an infinite sequence of 3-component links,   each of which is not zero on an invariant of us and is zero on $\mu_{123}$.   

\begin{figure}[H]
  \centering
  \includegraphics[width=5cm]{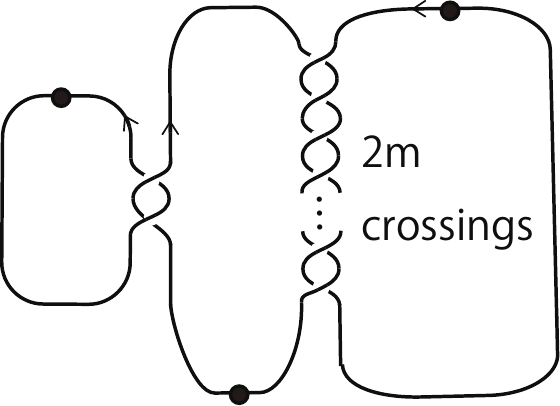}
  \caption{3-component link $L_{2m}^{+}$}\label{EgB}
\end{figure}

Let $L_{2m}^{+}$ be a 3-component link with base points as in Figure~\ref{EgB}.  Let $f_1$ be  an invariant $\F$ satisfying the assumption of (\ref{I}) of Theorem~\ref{main}.  
Then $\mu_{123}(L_{2m}^{+}) \equiv 0$ while $f_1 (L_{2m}^{+})  = m$ $( >0 )$. Likewise, For link $L_{2m}^{-}$ (Figure 2),  $\mu_{123}(L_{2m}^{-})  \equiv 0$ while $f_1(L_{2m}^{-}) = m$ $(<0)$.  \\

\begin{figure}[H]
  \centering
  \includegraphics[width=5cm]{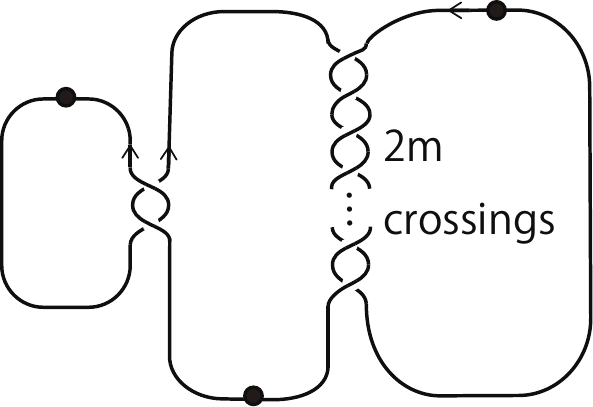}
  \caption{3-component link $L_{2m}^{-}$}\label{EgA}
\end{figure}

\section{Acknowlegements}
The work of NI was partially  supported by JSPS KAKENHI Grant number JP 20K03604.  

\bibliographystyle{plain}
\bibliography{Ref}

\begin{thebibliography}{1}

\bibitem{Deguchi1994}
Tetsuo Deguchi.
\newblock On numerical applications of the {V}assiliev invariants to
  computational problems in physics.
\newblock In {\em Proceedings of the {C}onference on {Q}uantum {T}opology
  ({M}anhattan, {KS}, 1993)}, pages 87--98. World Sci. Publ., River Edge, NJ,
  1994.

\bibitem{Ito2019}
Noboru Ito.
\newblock Space of chord diagrams on spherical curves.
\newblock {\em Internat. J. Math.}, 30(12):1950060, 25, 2019.

\bibitem{Ito2022}
Noboru Ito.
\newblock A triple coproduct of curves and knots.
\newblock {\em arXiv}, 2022.

\bibitem{ItoOyamaguchi2021}
Noboru Ito and Natsumi Oyamaguchi.
\newblock Gauss diagram formulas of {V}assiliev invariants of 2-bouquet graphs.
\newblock {\em Topology Appl.}, 290:Paper No. 107580, 13, 2021.

\bibitem{Ostlund2001PhD}
Olof-Petter Ostlund.
\newblock {\em Invariants of knot diagrams and diagrammatic knot invariants}.
\newblock ProQuest LLC, Ann Arbor, MI, 2001.
\newblock Thesis (Ph.D.)--Uppsala Universitet (Sweden).

\bibitem{Ostlund2004}
Olof-Petter \"{O}stlund.
\newblock A diagrammatic approach to link invariants of finite degree.
\newblock {\em Math. Scand.}, 94(2):295--319, 2004.

\bibitem{PolyakViro1994}
Michael Polyak and Oleg Viro.
\newblock Gauss diagram formulas for {V}assiliev invariants.
\newblock {\em Internat. Math. Res. Notices}, (11):445ff., approx. 8 pp.\,
  1994.

\end{thebibliography}
\end{document}